\documentclass[a4paper,11pt]{amsart}
\pdfoutput=1 
\usepackage{amsmath,amsthm,amssymb}
\usepackage[mathscr]{eucal}
 \usepackage{cite}
\usepackage{upgreek}
\usepackage[bookmarks=false]{hyperref}
\usepackage{enumerate}
\usepackage{amsmath}
\usepackage{mathrsfs}
\usepackage{blindtext}
\usepackage{scrextend}
\usepackage{enumitem}
\usepackage{bm} 
\addtokomafont{labelinglabel}{\sffamily}
\usepackage{color}
\usepackage[at]{easylist}
\usepackage{tikz}
\usepackage{svg}
\usepackage{graphicx}

\usepackage{comment}

\numberwithin{equation}{section}

\newtheorem{main}{Theorem}

\newtheorem*{thm*}{Theorem}

\newtheorem*{prob*}{Problem}

\newtheorem*{prop*}{Proposition}

\newtheorem*{cor*}{Corollary}

\theoremstyle{definition}

\newtheorem*{defn*}{Definition}

\newtheorem*{question*}{Question}
\newtheorem*{Pquestion*}{Popa's question}

\newtheorem*{conv*}{Convention}

\usepackage{tikz}
\usepackage{tikz-cd}

\begin{document}

\title[3-handle construction on II$_1$ factors]
{3-handle construction on II$_1$ factors}

\author[David Gao]{David Gao}
\address{Department of Mathematical Sciences, UCSD, 9500 Gilman Dr, La Jolla, CA 92092, USA}\email{weg002@ucsd.edu}\urladdr{https://sites.google.com/ucsd.edu/david-gao}

\author[Srivatsav Kunnawalkam Elayavalli]{Srivatsav Kunnawalkam Elayavalli}
\address{Department of Mathematical Sciences, UCSD, 9500 Gilman Dr, La Jolla, CA 92092, USA}\email{skunnawalkamelayaval}
\urladdr{https://sites.google.com/view/srivatsavke/home}

\author[Gregory Patchell]{Gregory Patchell}
\address{Department of Mathematical Sciences, UCSD, 9500 Gilman Dr, La Jolla, CA 92092, USA}\email{gpatchel@ucsd.edu}
\urladdr{https://sites.google.com/view/gpatchel/home}

\begin{abstract}
    We find a subtle modification to the construction in \cite{CIKE23}, which yields a drastic simplification of the proof of the existence of two non Gamma non elementarily equivalent II$_1$ factors.  
\end{abstract}

\maketitle

We introduce in this note a 3-handle construction on II$_1$ factors, which is directly inspired by and modifies the Chifan-Ioana-Kunnawalkam Elayavalli construction of an exotic non Gamma II$_1$ factor \cite{CIKE23} that gave the first concrete example of a pair of non Gamma, non elementarily equivalent II$_1$ factors. The novelty here is that the 3-handle construction makes the proofs drastically shorter, notably entirely removing two technical challenges of the \cite{CIKE23} approach: a deformation-rigidity argument to control relative commutants and ultrapower lifting theorems for pairs of independent unitaries. This new construction also yields sharper bounds on the commutation diameter \cite{elayavalli2023sequential} than previously known.

\subsection*{The construction} For a II$_1$ factor $N$ and $u,v\in \mathcal{U}(N)$, define $\Omega_1(N,u)= N*_{\{u\}''}\left(\{u\}''\mathbin{\bar{\otimes}} L(\mathbb{Z})\right)$  and denote the Haar unitary  generating the $L(\mathbb{Z})$ on the right as $\omega(u)$. Define  $\Omega_2(N,u,v)= N*_{\{u,v\}''}(\{u,v\}''\mathbin{\bar{\otimes}} 
 L(\mathbb{Z}))$. 
Let $\sigma=(\sigma_1,\sigma_2):\mathbb N\rightarrow\mathbb N\times\mathbb N$ be a bijection such that $\sigma_1(n)\leq n$, for every $n\in\mathbb N$. Assume that $M_1,\ldots,M_n$ have been constructed, for some $n\in\mathbb N$.
Let $\{(u_1^{n,k},u_2^{n,k})\}_{k\in\mathbb N}\subset \mathcal{U}(M_n)$ be a $\|\cdot\|_2$-dense sequence.
We define $$M_{n+1}:=\Omega_2(\Omega_1(\Omega_1(M_n,u_1^{\sigma(n)}),u_2^{\sigma(n)}), \omega(u_1^{\sigma(n)}), \omega(u_2^{\sigma(n)})).$$
Note that $M_{n+1}$ is well-defined since $\sigma_1(n)\leq n$ and thus $(u_1^{\sigma(n)},u_2^{\sigma(n)})\in\mathcal{U}(M_n)$. Denote by $M$ the inductive limit, where we start with $M_1= L(SL_3(\mathbb{Z}))$.  Note that since unitaries in ultrapowers lift into unitaries in strongly dense subsets, we have for all ultrafilters $U$ and any pair of unitaries $u,v\in M^U$, there exist diffuse unitaries $w_1,w_2,w_3\in M^U$ such that  $[u, w_1]= [w_1, w_2]= [w_2, w_3]= [w_3, v]=0$.

\subsection*{Showing that $M$ is a  non Gamma II$_1$ factor}
Observe that $L(SL_3(\mathbb{Z}))\subset M_n$ satisfies $L(SL_3(\mathbb{Z}))'\cap M_{n+1}= L(SL_3(\mathbb{Z}))'\cap M_{n}$. Indeed, by Theorem 1.1 of \cite{IPP08}, we have that if $L(SL_3(\mathbb{Z}))\not\prec_{M_{n+1}} \{\omega(u_1^{\sigma(n)}), \omega(u_2^{\sigma(n)})\}''$ then $L(SL_3(\mathbb{Z}))'\cap M_{n+1} \subset M_n$ (notice that the other two amalgams in the construction are abelian).   This is indeed the case: $L(SL_3(\mathbb{Z}))$ has Property (T) and the \emph{central point} of this note is that since $\omega(u_1^{\sigma(n)})$ and $\omega(u_2^{\sigma(n)})$ are by construction {freely independent} they generate a copy of the free group factor which has the Haagerup property (see \cite{ConnesJones} and Proposition 16.2.3 of \cite{anantharaman-popa}). Hence by weak spectral gap of Property (T) factors,  $L(SL_3(\mathbb{Z}))' \cap M^{\mathcal{U}}= \mathbb{C}$.

\subsection*{Proving non elementary equivalence} 

Let $a,b$ be the generator unitaries in $L(\mathbb F_2)$. For any ultrafilter $\mathcal{U}$ on a set $I$, it is known that there do not exist diffuse unitaries $w_1,w_2, w_3\in \mathcal{U}(L(\mathbb{F}_2)^\mathcal{U})$ such that $[a, w_1]= [w_1, w_2]= [w_2, w_3]= [w_3, b]=0$. There are now multiple proofs of the above statement. For instance, this follows immediately from repeatedly applying  the well known \emph{join lemma} in 1-bounded entropy theory (see Corollary 4.8 and Proposition A.12 of \cite{Hayes2018}, also \cite{JungSB},  see also  Fact 2.9 in \cite{elayavalli2023sequential}); alternatively, one can apply the Jekel-Kunnawalkam Elayavalli upgraded free independence theorem \cite{jekel2024upgradedfreeindependencephenomena}. Note that the free independence theorem of Houdayer-Ioana \cite{houdayer2023asymptotic} does not apply here as it can only handle one level of commutants.

The following theorem follows immediately by combining the  above with our 3-handle construction which by design attaches a handle of three sequentially commuting diffuse unitaries on every pair of unitaries in a countable strongly dense set.
\begin{main}
    For any ultrafilters $\mathcal{U},\mathcal{V}$ on sets $I, J$ respectively,   $M^\mathcal{U}\ncong L(\mathbb{F}_2)^\mathcal{V}$. Moreover $3\leq \mathfrak{D}(M)\leq 4$ where $\mathfrak{D}(M)$ denotes the commutation diameter introduced in \cite{elayavalli2023sequential}. 
\end{main}

\bibliographystyle{amsalpha}
\bibliography{inneramen}

\providecommand{\bysame}{\leavevmode\hbox to3em{\hrulefill}\thinspace}
\providecommand{\MR}{\relax\ifhmode\unskip\space\fi MR }
\providecommand{\MRhref}[2]{%
  \href{http://www.ams.org/mathscinet-getitem?mr=#1}{#2}
}
\providecommand{\href}[2]{#2}
\begin{thebibliography}{CIKE23}

\bibitem[AP16]{anantharaman-popa}
Claire Anantharaman and Sorin Popa, \emph{An introduction to {$II_1$} factors}, book in progress (2016).

\bibitem[CIKE23]{CIKE23}
Ionu\c{t} Chifan, Adrian Ioana, and Srivatsav Kunnawalkam~Elayavalli, \emph{An exotic {$\rm II_1$} factor without property gamma}, Geom. Funct. Anal. \textbf{33} (2023), no.~5, 1243--1265. \MR{4646408}

\bibitem[CJ85]{ConnesJones}
Alain Connes and Vaughan Jones, \emph{Property {$T$} for von {N}eumann algebras}, Bull. London Math. Soc. \textbf{17} (1985), no.~1, 57--62. \MR{766450}

\bibitem[Hay18]{Hayes2018}
Ben Hayes, \emph{1-bounded entropy and regularity problems in von {N}eumann algebras}, Int. Math. Res. Not. IMRN (2018), no.~1, 57--137. \MR{3801429}

\bibitem[HI23]{houdayer2023asymptotic}
Cyril Houdayer and Adrian Ioana, \emph{Asymptotic freeness in tracial ultraproducts}, 2023.

\bibitem[IPP08]{IPP08}
Adrian Ioana, Jesse Peterson, and Sorin Popa, \emph{Amalgamated free products of weakly rigid factors and calculation of their symmetry groups}, Acta Math. \textbf{200} (2008), no.~1, 85--153. \MR{2386109}

\bibitem[JKE24]{jekel2024upgradedfreeindependencephenomena}
David Jekel and Srivatsav Kunnawalkam~Elayavalli, \emph{Upgraded free independence phenomena for random unitaries}, 2024.

\bibitem[Jun07]{JungSB}
Kenley Jung, \emph{Strongly 1-bounded von {N}eumann algebras}, Geom. Funct. Anal. \textbf{17} (2007), no.~4, 1180--1200. \MR{2373014}

\bibitem[KEP25]{elayavalli2023sequential}
Srivatsav Kunnawalkam~Elayavalli and Gregory Patchell, \emph{Sequential commutation in tracial von {N}eumann algebras}, J. Funct. Anal. \textbf{288} (2025), no.~4, Paper No. 110719, 28. \MR{4832103}

\end{thebibliography}

\end{document}